\documentclass[amssymb,amstex,amsmath]{amsart}

\usepackage{amsmath}
\usepackage{amssymb}
\usepackage{epsfig}

\newtheorem{thm}{Theorem}

\newenvironment{pf}{\noindent \rm Proof.}{$\hspace{8mm}\square$}

\begin{document}

\title[Embedded minimal and cmc annulus touching spheres]
{Embedded minimal and constant mean curvature annulus touching
spheres}

\author[S.-H. Park]{Sung-Ho Park}
\address{Korea Institute for Advanced Study, 207-43 Cheongryangri
2-dong, Dongdaemun-gu, Seoul 130-722,
 Korea}
\email{shubuti@kias.re.kr}

\begin{abstract}
We show that a compact embedded constant mean curvature annulus in
$\mathbb R^3$ tangent to two spheres of same radius along its
boundary curves and having non-vanishing Gaussian curvature is
part of a Delaunay surface. In special, if the annulus is minimal,
then the annulus is part of a catenoid. Secondly we show that a
compact embedded constant mean curvature annulus with negative
(respectively, positive) Gaussian curvature meeting a sphere
tangentially and a plane in constant contact angle $\ge\pi/2$
(respectively, $\le\pi/2$) is part a Delaunay surface. In special,
if the annulus is minimal and the contact angle is $\ge\pi/2$,
then it is part of a catenoid.\\

{\it 2010 Mathematics Subject Classification}: 53A10.\\
\end{abstract}

\maketitle Catenoid is the only nonplanar minimal surface of
rotation in $\mathbb R^3$ \cite{bo}. Rotational surfaces of
constant mean curvature in $\mathbb R^3$ are called the Delaunay
surfaces: cylinders, spheres, unduloids and nodoids. Therefore
catenoid and Delaunay surface meet every plane, which is
perpendicular to the axis of rotation, in constant contact angle.
Conversely, if a compact embedded minimal or constant mean
curvature surface meets two parallel planes in constant contact
angles, then the surface is part of a catenoid or part of a
Delaunay surface. This can be proved by using the Alexandrov's
moving plane argument \cite{hopf} to planes perpendicular to the
parallel planes. A compact immersed minimal annulus meeting two
parallel planes in constant contact angles is also part of a
catenoid. This result is not true for constant mean curvature
surfaces: Wente had constructed examples of immersed constant mean
curvature annuli in a slab or in a ball meeting the boundary
planes or the boundary sphere perpendicularly \cite{W}. Compared
to the above first case, we may ask whether a compact minimal
annulus or a compact embedded constant mean curvature annulus
meeting two spheres in constant contact angles is part of a
catenoid or part of a plane. In \cite{pp}, it is shown that if a
compact embedded minimal annulus meets two concentric spheres
perpendicularly then the minimal annulus is part of a plane.

In this paper, we show that a compact embedded constant mean
curvature annulus $\mathcal{A}$ in $\mathbb R^3$ meeting two
spheres $S_1$ and $S_2$ of same radius $\rho$ tangentially and
having non-vanishing Gaussian curvature $K$ is part of a Delaunay
surface. More precisely, depending on the values of $K$ and the
mean curvature $H$ we have three cases: i) $K<0$ and $H>-1/\rho$,
in which case $\mathcal{A}$ is part of a unduloid if $H<0$, part
of a catenoid if $H=0$ and part of a nodoid if $H>0$, ii) $K>0$
and $-1/\rho<H<-1/2\rho$, in which case $\mathcal{A}$ is part of a
unduloid, and iii) $K>0$ and $H<-1/\rho$, in which case
$\mathcal{A}$ is part of a nodoid. In the first two cases,
$\mathcal{A}$ stays outside of the balls $B_1$ and $B_2$ bounded
by $S_1$ and $S_2$. If iii) holds, then $\mathcal{A} \subset
B_1\cap B_2$.

We also show that a compact embedded constant mean curvature
annulus $\mathcal{B}$ in $\mathbb R^3$ with negative
(respectively, positive) Gaussian curvature meeting a unit sphere
tangentially and a plane in constant contact angle $\ge\pi/2$
(respectively, $\le\pi/2$) is part of a Delaunay surface. In
special, a compact embedded minimal annulus in $\mathbb R^3$
meeting a sphere tangentially and a plane in constant contact
angle $\ge\pi/2$ is part of a catenoid.

To prove Theorem 1 and 2, we use the $-\rho$-parallel surface
$\tilde{\mathcal{A}}$ of $\mathcal{A}$ (respectively,
$\tilde{\mathcal{B}}$ of $\mathcal{B}$): the parallel surface of
$\mathcal{A}$ (respectively, of $\mathcal{B}$) with distance
$\rho$ in the direction to the centers of the spheres. We use the
Alexandrov's moving plane argument \cite{al2}, \cite{hopf} to
prove that $\tilde{\mathcal{A}}$ and $\tilde{\mathcal{B}}$ are
rotational. Since $\tilde{\mathcal{A}}$ and $\tilde{\mathcal{B}}$
are the parallel surfaces of $\mathcal{A}$ and $\mathcal{B}$
respectively, $\mathcal{A}$ and $\mathcal{B}$ are also rotational
and, hence, are part of a Delaunay surface or part of a catenoid.

\section{constant mean curvature annulus meeting spheres
tangentially} In the following, we may assume that the spheres
have radius $1$. Let $\mathcal{A}$ be a compact embedded annulus
with constant mean curvature $H$ and meeting two unit spheres
$S_1$ and $S_2$ tangentially along the boundary curves $\gamma_1$
and $\gamma_2$. We fix the unit normal $N$ of $\mathcal{A}$ to
point away from the centers of the spheres. Let $Y: A(1, R)
\rightarrow \mathbb R^3$ be a conformal parametrization of
$\mathcal{A}$ from an annulus $A(1,R)= \{(x,y)\in \mathbb R^2:
1\le \sqrt{x^2 +y^2} \le R\}$. We define $X$ by $X= Y \circ \exp$
on the strip $B=\{(u,v)\in \mathbb R^2: 0\le u\le \log R\}$. Then
$X$ is periodic with period $2\pi$. Let $z=u+iv$ and $\lambda^2 :=
|X_u| ^2 = |X_v|^2$.

Let $h_{ij}$, $i,j=1,2$, be the coefficients of the second
fundamental form of $X$ with respect to $N$. Note that the Hopf
differential  $\phi(z)dz^2=(h_{11}- h_{22} -2i h_{12})dz^2$ is
holomorphic for constant mean curvature surfaces \cite{hopf}. The
theorem of Joachimstahl \cite{dc} says that $\gamma_1$ and
$\gamma_2$ are curvature lines of $\mathcal{A}$. Hence
$h_{12}\equiv 0$ on $u=0$ and $u= \log R$. Since $h_{12}$ is
harmonic and periodic, we have $h_{12}\equiv 0$ on $B$. This
implies that $z$ is a conformal curvature coordinate and $h_{11}-
h_{22}=constant$ \cite{mc}. Let $c=h_{11}-h_{22}$. If
$\mathcal{A}$ is minimal, then we have $K<0$ and $c=2h_{11}>0$ by
the choice of $N$. When $H=-1$, $\mathcal{A}$ is part of the unit
sphere $S_1=S_2$ by the boundary comparison principle for mean
curvature operator \cite{gt}. We assume that $H\not=-1$ in the
following. The principal curvatures of $\mathcal{A}$ are
\begin{equation}\label{pcurv}
\kappa_1= H+ {c \over 2\lambda^2} \mbox{ and } \kappa_2= H- {c
\over 2\lambda^2}.
\end{equation}

We parameterize $\gamma_1$ and $\gamma_2$ by $\gamma_1(v)= X(0,v)$
and $\gamma_2(v)= X(\log R,v)$ for $v\in [0,2\pi)$. In the
following, we assume that $\mathcal{A}$ has non-vanishing Gaussian
curvature.\\

\noindent{\bf Lemma 1.} {\it Each $\gamma_i(v)$, $i=1,2$, has
constant speed $\sqrt{c/ 2(1+H)}$ and $\kappa_2$ is $-1$ on
$\gamma_1$ and $\gamma_2$. As spherical curves, $\gamma_1$ and
$\gamma_2$ are convex. On
$\mathcal{A}\setminus\partial\mathcal{A}$, we have $\lambda^2<{c/
2(1+H)}$ when $K<0$ and $\lambda^2>{c/2(1+H)}$ when $K>0$.} \\

\begin{pf}
The curvature vector of $\gamma_1(v)$ is
\begin{eqnarray}\label{curvature}
\vec{\kappa} &=& {1\over |X_v|}{d \over dv}\left( X_{v}\over |X_v|
\right) = {1\over |X_v|^2}X_{vv} - { X_{v}\over |X_v|^4}(X_v\cdot
X_{vv}) \\
&=& {1 \over \lambda^2} \left( -{\lambda_u \over \lambda} X_u +
h_{22}N\right). \nonumber
\end{eqnarray}
Let the center of $S_1$ be the origin of $\mathbb R^3$. Since
$\mathcal{A}$ is tangential to $S_1$ along $\gamma_1$, we have
$N(0,v)=X(0,v)=\gamma_1(v)$ on $\gamma_1$. Since $\gamma_1$ is on
the unit sphere $S_1$, the curvature vector $\vec{\kappa}$ of
$\gamma_1$ satisfies $(\vec{\kappa}\cdot \gamma_1)(v)=-1$. Hence
we have $\kappa_2={h_{22}\over \lambda^2} =-1$ on $\gamma_1$.
Since $\lambda^2 = |{\gamma_1}_v|^2$ on $\gamma_1$, we have
$|{\gamma_1}_v|=\sqrt{c/ 2(1+H)}$ from \eqref{pcurv}. By choosing
the center of $S_2$ as the origin of $\mathbb R^3$, we get the
results for $\gamma_2$.

The Gaussian curvature $K$ satisfies
\[
\Delta \log \lambda = - K \lambda^2,
\] where $\Delta = {\partial ^2 \over \partial u^2} +{\partial ^2
\over \partial v^2}$. We can rewrite this equation as
\begin{equation}\label{gc}
\lambda \Delta \lambda = |\nabla \lambda|^2 -K \lambda^4.
\end{equation}
Since $\lambda_v(0,v)=0$ and $\lambda_v(\log R,v) =0$ and
$K\not=0$, $\lambda$ does not have interior maximum when $K<0$,
and does not have interior minimum when $K>0$. Since
$\lambda^2={c/2(1+H)}$ on $\gamma_1$ and $\gamma_2$, it follows
that $\lambda^2<{c/2(1+H)}$ on
$\mathcal{A}\setminus\partial\mathcal{A}$ when $K<0$ and
$\lambda^2>{c/2(1+H)}$ when $K>0$. Moreover we have
$\lambda_u\le0$ on $u=0$ and $\lambda_u\ge 0$ on $u=\log R$ when
$K<0$ and $\lambda_u\ge0$ on $u=0$ and $\lambda_u\le 0$ on $u=\log
R$ when $K>0$. Since ${X_u\over |X_u|}\in TS_i$ is perpendicular
to $\gamma_i$, the geodesic curvature of $\gamma_i$ as a spherical
curve is
$\vec{\kappa}\cdot{X_u\over|X_u|}=-{\lambda_u\over\lambda^2}$.
Hence $\gamma_1$ and $\gamma_2$ are convex as spherical curves.
\end{pf}\\

\noindent{\bf Remark 1.} If $\lambda^2\equiv {c/2(1+H)}$ on
$\mathcal{A}$, then $K\equiv0$ and $\mathcal{A}$ is part of a
cylinder.

\section{$-1$-parallel surface}
The $-1$-parallel surface $\tilde{\mathcal{A}}$ of $\mathcal{A}$
is defined by
\[
\tilde{X} = X - N.
\]
The image of $\gamma_1$ (respectively, of $\gamma_2$) in
$\tilde{\mathcal{A}}$ is a point corresponding to the center of
$S_1$ (respectively, of $S_2$). We denote the centers of $S_1$ and
$S_2$ by $O$ and $O_2$ for simplicity. We fix the unit normal
$\tilde{N}$ of $\tilde{\mathcal{A}}$ to be $N$. Since $z=u+iv$ is
a curvature coordinate of $X$, we have
\begin{equation}\label{derivative}
\tilde{X}_u = \left(1+ {h_{11} \over \lambda^2}\right) X_u \mbox{
and } \tilde{X}_v = \left(1+ {h_{22} \over \lambda^2}\right) X_v.
\end{equation}
Since $\kappa_2=-1$ on $\gamma_i$ (Lemma 1), $\tilde{X}$ is
singular for $u=0$ and $u=\log R$. By Lemma 1, we have
$\lambda^2\not= {c/2(1+H)}$ on
$\mathcal{A}\setminus\partial\mathcal{A}$, which implies that
$1+\kappa_2\not=0$ on $\mathcal{A}\setminus\partial\mathcal{A}$.
When $K<0$, we have $\kappa_1> 0$ on
$\mathcal{A}\setminus\partial\mathcal{A}$. Hence $\tilde{X}$ is
regular for $0<u<\log R$ and we have $H>-1$.

Now suppose that $K>0$. Since $\kappa_2=-1$ on $\gamma_i$ (Lemma
1), we have $\kappa_1<0$ and $H<-1/2$. We consider two cases
separately: $H<-1$ and $-1<H<-1/2$. If $H<-1$, then $c<0$
from $\lambda^2 ={c/2(1+H)}>0$ on $\gamma_i$. Hence we have
$\kappa_1<-1$, which implies that $\tilde{X}$ is regular for
$0<u<\log R$. If $-1<H<-1/2$, then we must have $c>0$. This
implies that $1+\kappa_1\not=0$. Otherwise we have
$0<2\lambda^2(1+H) =-c$, which contradicts $c>0$. Hence
$\tilde{X}$ is regular for $0<u<\log R$.\\

\noindent{\bf Remark 2.} When $K<0$ or $K>0$ and $-1<H<-1/2$,
$\mathcal{A}$ stays outside of the balls $B_1$ and $B_2$ bounded
by $S_1$ and $S_2$. If $K>0$ and $H<-1$, then
$\mathcal{A}\subset B_1\cap B_2$.\\

\noindent{\bf Lemma 2.} {\it The mean curvature $\tilde{H}$ and
the Gaussian curvature $\tilde{K}$ of $\tilde{\mathcal{A}}$
satisfies $(1+H)\tilde{K}=(1+2H)\tilde{H} -H$. On
$\tilde{\mathcal{A}}\setminus\{O,O_2\}$, we have

i) if $K<0$ and $H>-1$, then $\tilde{\kappa}_1>0$,
$\tilde{\kappa}_2>1$ and $\tilde{H}>1$,

ii) if $K>0$ and $-1<H<-1/2$, then $0<c/2\lambda^2(1+H) <\min\{1,
-H/(1+H)\}$, $\tilde{\kappa}_1<0$, $\tilde{\kappa}_2<H/(1+H)$ and
$\tilde{H}<H/(1+H)$, and

iii) if $K>0$ and $H<-1$, then $0<c/2\lambda^2(1+H) <1$,
$\tilde{\kappa}_1>(1+2H)/2(1+H)$, $\tilde{\kappa}_2>H/(1+H)$ and
$\tilde{H}>H/(1+H)$.}\\

\begin{pf}
Since
\[
\tilde{h}_{12}=N \cdot \tilde{X}_{uv} = \left(1+ {h_{11} \over
\lambda^2}\right) (N\cdot X_{uv}) =0,
\]
$(u,v)$ is a curvature coordinate (not conformal) for
$\tilde{\mathcal{A}}$ except for $O$ and $O_2$.
We have
\begin{eqnarray*}
\tilde{h}_{11}=N \cdot \tilde{X}_{uu} = \left(1+ {h_{11} \over
\lambda^2}\right)h_{11},\\ \tilde{h}_{22}= N \cdot \tilde{X}_{vv}
=\left(1+ {h_{22} \over \lambda^2}\right)h_{22}.
\end{eqnarray*}
The principal curvatures of $\tilde{\mathcal{A}}$ are
\begin{eqnarray*}
\tilde{\kappa}_1=  {\kappa_{1} \over 1+ \kappa_1} =
{{H/(1+H)}+\left(c/2\lambda^2(1+H)\right)\over
{1+\left(c/2\lambda^2(1+H)\right)}},\\
\tilde{\kappa}_2=  {\kappa_{2} \over 1+ \kappa_2}=
{{{H/(1+H)}-\left(c/2\lambda^2(1+H)\right)}\over
{1-\left(c/2\lambda^2(1+H)\right)}}.
\end{eqnarray*}
From $\kappa_1+\kappa_2 =2H$, we have
\[
H= {\tilde{H}- \tilde{K} \over 1-2\tilde{H} -\tilde{K}}\mbox{
or } (1+H)\tilde{K}=(1+2H)\tilde{H} -H.
\]
It is straightforward to see that \[
\tilde{H}= {{{H/(1+H)}-\left(c/2\lambda^2(1+H)\right)^2}\over
{1-\left(c/2\lambda^2(1+H)\right)^2}}.
\]

Note that $\kappa_2<0$ on $\mathcal{A}$. First suppose that $K<0$.
Then we have $\kappa_1>0$, which implies that $\tilde{\kappa}_1=
\kappa_1/(1+\kappa_1)>0$. Since $c/2\lambda^2(1+H)>1$ by Lemma 1,
we have $\tilde{\kappa}_2>1$ and $\tilde{H}>1$.

When $K>0$, we have $\kappa_1= H+c/2\lambda^2<0$. If $-1<H<-1/2$,
then we have $c>0$ because $\lambda^2=c/2(1+H)>0$ on $\gamma_i$.
It follows that $c/2\lambda^2(1+H)<-H/(1+H)$. By Lemma 1, we also
have $c/2\lambda^2(1+H)<1$. Therefore we have
$0<c/2\lambda^2(1+H)<\min\{1,-H/(1+H)\}$. It is straightforward to
see that $\tilde{\kappa}_1<0$ and $\tilde{\kappa}_2<H/(1+H)<0$ and
$\tilde{H}<H/(1+H)<0$.

When $K>0$ and $H<-1$, we have $c<0$ and $0<c/2\lambda^2(1+H)<1$.
It is straightforward to see that $\tilde{\kappa}_1>(1+2H)/(1+H)$,
$\tilde{\kappa}_2>H/(1+H)$ and $\tilde{H}>H/(1+H)$.
\end{pf}\\

This lemma says that $\tilde{\mathcal{A}}$ is a linear Weingarten
surface with two singular points $O$ and ${O}_2$ and is positively
curved outside $O$ and $O_2$.\\

\noindent{\bf Lemma 3.} {\it $\tilde{\mathcal{A}}$ is embedded.}\\

\begin{pf} Let $\nu(v)={X_u\over|X_u|}(0,v)$.
Note that $\nu$ is a closed curve in the unit sphere $S_1$. We
claim that $\nu$ is {\it convex as a spherical curve}. Otherwise,
there is a great circle $\eta$ intersecting the image of $\nu$ at
no less than $3$ points $\nu(v_1),\ldots,\nu(v_n)$. ($\nu$ may map
an interval $(v_a, v_b)\subset[0, 2\pi)$ into a single point. We
choose $v_i$'s in such a way that $\nu$ maps no two $v_i$'s to the
same point.) Each $\nu(v_i)$ determines a great circle $\mathbb
S^1_{v_i} \subset S_1$ contained in the plane perpendicular to
$\nu(v_i)$. At each $\gamma_1(v_i)$, $\gamma_1$ is tangent to
$\mathbb S^1_{v_i}$. Since $\eta$ and $\mathbb S^1_{v_i}$ are
perpendicular, $\gamma_1$ cannot be convex when $n\ge3$. Hence
$\nu$ intersect every geodesic of $S_1$ at no more than two
points. This shows that $\nu$ is convex as a spherical curve.
Similarly, ${X_u\over|X_u|}(\log R,v)$ is also convex as a
spherical curve.

Since $\tilde{\mathcal{A}}$ is a parallel surface of
$\mathcal{A}$, the tangent cone $Tan(O,\tilde{\mathcal{A}})$ of
$\tilde{\mathcal{A}}$ at $O$ is the cone formed by rays from $O$
through $\nu$. Since $\nu$ is a convex spherical curve,
$Tan(O,\tilde{\mathcal{A}})$ is convex. This shows that a small
neighborhood of $O$ in $\tilde{\mathcal{A}}$ is embedded and
nonnegatively curved as a metric space \cite{al1}. Similarly,
there is a neighborhood of $O_2$ in $\tilde{\mathcal{A}}$ which is
embedded and nonnegatively curved as a metric space.

Hadamard showed that a closed surface in $\mathbb R^3$ with
strictly positive Gaussian curvature is the boundary of a convex
body \cite{hopf}. In particular, $S$ is embedded. Alexandrov
generalized Hadamard's theorem to nonnegatively curved metric
spaces \cite{al1}. Since $\tilde{\mathcal{A}}$ is a nonnegatively
curved closed metric space, $\tilde{\mathcal{A}}$ is embedded.
\end{pf}\\

\noindent{\bf Remark 2.} We have $\nu_v= {\lambda_u \over
\lambda^2} X_v$. At points where $\lambda_u \not=0$, the curvature
vector of $\nu$ is
\[
\vec{\kappa}_\nu = {1 \over \lambda_u} \left( -{\lambda_u \over
\lambda} X_u + h_{22}N\right).
\]
The geodesic curvature of $\nu$ as a spherical curve
$\vec{\kappa}_\nu\cdot N ={h_{22}\over \lambda_u}$.\\

\section{Main results}

We use the Alexandrov's moving plane argument \cite{al2},
\cite{hopf} to prove the theorems.

\begin{thm}
A compact embedded constant mean curvature annulus $\mathcal{A}$
with non-vanishing Gaussian curvature and meeting two spheres
$S_1$ and $S_2$ of same radius tangentially is part of a Delaunay
surface. In special, if $\mathcal{A}$ is minimal, then
$\mathcal{A}$ is part of a catenoid.
\end{thm}
\begin{pf}
We suppose that the radius of $S_1$ and $S_2$ is $1$. By Lemma 2
and Lemma 3, $\tilde{\mathcal{A}}$ is a compact embedded surface
with two singular points $O$ and $O_2$ and satisfying
$(1+H)\tilde{K}=(1+2H)\tilde{H} -H$ at regular points. A small
neighborhood of a regular point of $\tilde{\mathcal{A}}$ can be
represented as the graph of a function $f(x,y)$ satisfying
\begin{eqnarray}\label{pde}
&& 2(1+H)(f_{xx}f_{yy}-f_{xy}^2) +2H(1+f_x^2+f_y^2)^2
\\ &&= (1+2H)\left((1+f_y^2)f_{xx}-2f_xf_y
f_{xy}+(1+f_x^2)f_{yy}\right)(1+f_x^2+f_y^2)^{1\over2}. \nonumber
\end{eqnarray}
This equation can be rewritten as
\begin{equation}\label{pde2}
\det\left(2(1+H)D^2f + A(Df)\right)=W^4,
\end{equation}
where $A(Df)=-(1+2H)\left(\begin{array}{c}(1+f_x^2)W\\
f_x f_y W\end{array}\begin{array}{c}f_x f_y W\\ (1+f_y^2)W
\end{array}\right)$ and $W=\sqrt{1+f_x^2 +f_y^2}$. The equation
\eqref{pde2} is elliptic with respect to $f$ if $2(1+H)D^2f+A(Df)$
is positive definite. Since
$\det\left(2(1+H)D^2f+A(Df)\right)=W^4>0$, \eqref{pde2} is
elliptic if
\begin{equation}\label{ec1}
\mbox{Tr}\left(2(1+H)D^2f+A(Df)\right)= 2(1+H)\Delta f
-(1+2H)(2+ f_x^2 +f_y^2)W
\end{equation}
is strictly positive.

First we consider the case $K<0$. Since $\tilde{H}>1$ by Lemma 2,
we have
\begin{equation}\label{ec2}
\Delta f + f_y^2f_{xx} -2f_x f_y f_{xy} + f_x^2f_{yy}>2W^{3/2},
\end{equation}
for $f$ representing $\tilde{\mathcal{A}}$. We may assume that $f$
is defined on $B(0,\epsilon) \subset T_p\tilde{\mathcal{A}}$ so
that $\nabla f(0)=\vec{0}$ and $D^2f$ is diagonal. For
sufficiently small $\epsilon=\epsilon(p)$, \eqref{ec2} implies
that \eqref{ec1} is strictly positive. Hence \eqref{pde2} is
elliptic with respect to $f$ representing $\tilde{\mathcal{A}}$.

When $-1<H<-1/2$, \eqref{ec1} is automatically satisfied.

Now we consider the case $K>0$ and $H<-1$. Since
$\tilde{H}>H/(1+H)$ by Lemma 2, we have
\begin{equation}\label{ec3}
\Delta
f+f_y^2f_{xx}-2f_xf_yf_{xy}+f_x^2f_{yy}>{2H\over1+H}W^{3/2}.
\end{equation}
Assuming that $f$ is defined on $B(0,\epsilon)\subset
T_p\tilde{\mathcal{A}}$ with $\nabla f(0)=\vec{0}$ and $D^2f$ is
diagonal, \eqref{ec3} implies that $\Delta f-{1+2H\over2(1+H)}(2+
f_x^2 +f_y^2)W$ is strictly positive for sufficiently small
$\epsilon$. Then $\det\left(-2(1+H)D^2f-A(Df)\right)=W^4$ is
elliptic for $f$ representing $\tilde{A}$. The ellipticity of
\eqref{pde2} for $f$ representing $\tilde{A}$ enables us to use
the maximum principle and the boundary point lemma \cite{gt}.

Since $\tilde{\mathcal{A}}$ is convex and embedded, we can use
Alexandrov's moving plane argument \cite{al2}, \cite{hopf} to show
that $\tilde{\mathcal{A}}$ is rotational as follows. Let
$\Pi_\theta$ be the plane containing the line segment $\overline{O
O}_2\subset \mathbb R^3$ and making angle $\theta$ with a fixed
vector $\vec{E}$ which is perpendicular to $\overline{OO}_2$. Fix
a positive constant $L$ such that each plane $\Pi_\theta^{L}$,
which is parallel to $\Pi_\theta$ with distance $L$ from
$\Pi_\theta$, does not meet $\tilde{\mathcal{A}}$ for all
$\theta$. Let $\Pi_\theta^{l}$ be the plane between
$\Pi_\theta^{L}$ and $\Pi_\theta$ with distance $l$ from
$\Pi_\theta$. When $\Pi_\theta^{l}$ intersects
$\tilde{\mathcal{A}}$, we reflect the $\Pi_\theta^{L}$ side part
of $\tilde{\mathcal{A}}$ about $\Pi_\theta^{l}$. Let us denote
this reflected surface $\tilde{\mathcal{A}}_{l,\theta}^{ref}$. As
we decrease $l$ from $L$, there might be the first $l_\theta\ge 0$
for which $\tilde{\mathcal{A}}_{l_\theta,\theta}^{ref}$ is tangent
to $\tilde{\mathcal{A}}$ at an interior point or at a boundary
point of $\partial\tilde{\mathcal{A}}_{l_\theta,\theta}^{ref}$. We
call this point as the {\it first touch point}. If there is no
nonnegative $l$ with the first touch point, we repeat the process
for $\Pi_{\theta+ \pi}^{L}$ to find $l_{\theta+ \pi}$, which must
be positive. At the first touch point, we apply the comparison
principles for \eqref{pde} to see that the part of
$\tilde{\mathcal{A}}$ in the $\Pi_\theta$ side and
$\tilde{\mathcal{A}}_{l_\theta,\theta}^{ref}$ are identical and,
hence, $l_\theta =0$. This implies that $\Pi_\theta$ is a symmetry
plane for $\tilde{\mathcal{A}}$. Since $\theta$ can be chosen
arbitrarily, $\tilde{\mathcal{A}}$ should be rotational and,
hence, $\mathcal{A}$ is also rotational. Since the Delaunay
surfaces and the catenoid are the only nonplanar rotational
minimal and constant mean curvature surfaces, $\mathcal{A}$ is
part of a a Delaunay surface or part of a catenoid.
\end{pf}\\

We used the embeddedness of $\mathcal{A}$ in proving that
$\tilde{\mathcal{A}}$ is embedded. Whether there is a non-embedded
minimal or constant mean curvature annulus meeting two unit
spheres tangentially is an interesting question. Moreover we raise
the following questions.

1. Is a compact immersed minimal annulus or a compact embedded
minimal or constant mean curvature surface meeting a sphere
perpendicularly or in constant contact angles part of a catenoid
or part of a Delaunay surface? Nitsche showed that an immersed
disk type minimal or constant mean curvature surface meeting a
sphere in constant contact angle is either a flat disk or a
spherical cap \cite{N}.

2. Is a compact immersed minimal annulus or a compact embedded
minimal or constant mean curvature surface meeting two spheres in
constant contact angles part of a catenoid or a plane or part of a
Delaunay surface?

3. Is a compact immersed minimal or constant mean curvature
annulus or a compact embedded minimal or constant mean curvature
surface meeting a sphere and a plane in constant contact angles
part of a catenoid or part of a Delaunay surface? We give an
affirmative answer to this problem in a special case in the
following.

\begin{thm}
A compact embedded constant mean curvature annulus $\mathcal{B}$
with negative (respectively, positive) Gaussian curvature meeting
a sphere tangentially and a plane in constant contact angle $\ge
\pi/2$ (respectively, $\le \pi/2$) is part of a Delaunay surface.
In special, if $\mathcal{B}$ is minimal and the constant contact
angle is $\ge \pi/2$ then $\mathcal{B}$ is part of a catenoid.
\end{thm}

The angle is measured between the outward conormal of
$\mathcal{B}$ and the outward conormal of the bounded domain in
$\Pi$ bounded by the boundary curve. Since the proof of this
theorem is similar to that
of Theorem 1, we omit some details which was previously proved.\\

\begin{pf} Let us denote the sphere by $S_2$ and the plane by
$\Pi$. We may assume that the radius of $S_2$ is $1$. Let $\alpha$
be the constant contact angle between $\mathcal{B}$ and $\Pi$. If
$\alpha =\pi/2$, then we can reflect $\mathcal{B}$ about $\Pi$ to
get a constant mean curvature annulus meeting two unit spheres
tangentially. Hence $\mathcal{B}$ is part of a catenoid or a
Delaunay surface by Theorem 1.

In the following, we assume that $\alpha \not=\pi/2$. As in the
case for $\mathcal{A}$ in $\S 1$, there is a conformal
parametrization $X$ of $\mathcal{B}$ from a strip $\{(u,v)\in
\mathbb R^2:0\le u\le\log R\}$ for which $z=u+iv$ is a curvature
coordinate. We fix the normal $N$ of $\mathcal{B}$ to point away
from the center of $S_2$. Let $c_1(v)=X(0,v)$ be on $\Pi$ and
$c_2(v)=X(\log R,v)$ be on $S_2$ with $\partial X_3/\partial u>0$
along $c_1$. As in Lemma 1, $c_2$ has constant speed
$\sqrt{c/2(1+H)}$ and $\kappa_2=-1$ along $c_2$. Since $K\not=0$
on $\mathcal{B}$ and $z=u+iv$ is a curvature coordinate, we have
$\kappa_2<0$ on $c_1$. The curvature of $c_1$ is
$|\vec{\kappa}|=-\kappa_2/\sin\alpha>0$, which shows that $c_1$ is
locally convex. Since $c_1$ is a Jordan curve, it is convex.

First, we assume that $K<0$ and $\alpha >\pi/2$. Since
${\vec{\kappa}\over |\vec{\kappa}|}\cdot {X_u\over|X_u|}= \cos
\alpha <0$ on $c_1$, it follows from \eqref{curvature} that
$\lambda_u>0$ on $c_1$. Since $\lambda_v(\log R,v)=0$ (cf. Lemma
1), it follows from \eqref{gc} that $\lambda_u\ge0$ on $c_2$.
Otherwise, $\lambda$ will have an interior maximum, which
contradicts \eqref{gc}. Hence we have $\lambda^2<c/2(1+H)$ on
$\mathcal{B}\setminus c_2$. Note that $\kappa_1>0$ and
$\kappa_2<0$ in $\mathcal{B}$. From $\lambda_u\le0$ on $c_2$, we
see that $c_2$ is convex as a spherical curve (cf. Lemma 1).
Arguing as in the proof of Lemma 3, we see that
${X_u\over|X_u|}(\log R,v)$ is also convex as a spherical curve.

When $K>0$ and $\alpha <\pi/2$, we have ${\vec{\kappa}\over
|\vec{\kappa}|}\cdot {X_u\over|X_u|}= \cos \alpha >0$ on $c_1$.
Hence $\lambda_u<0$ on $c_1$. Since $\lambda_v(\log R,v)=0$, it
follows from \eqref{gc} that $\lambda$ does not have interior
minimum. Then we have $\lambda_u\le0$ on $c_2$ and
$\lambda^2>c/2(1+H)$ on $\mathcal{B} \setminus c_2$. Note that
$\kappa_1<0$ and $\kappa_2<0$ in $\mathcal{B}$. From
$\lambda_u\le0$ on $c_2$, it follows that $c_2$ is convex as a
spherical curve. Moreover ${X_u\over|X_u|}(\log R,v)$ is convex as
a spherical curve (cf. Lemma 3).

Let $\tilde{\mathcal{B}}$ be the $-1$-parallel surface of
$\mathcal{B}$. As in $\S 2$, we can show that
$\tilde{\mathcal{B}}$ is regular except for $O_2$: the image of
$c_2$, and $H>-1$ when $K<0$ and $H<-1/2$ when $K>0$. As in Lemma
2, we see that mean curvature $\tilde{H}$ and the Gaussian
curvature $\tilde{K}$ of $\tilde{\mathcal{B}}$ satisfies
$(1+H)\tilde{K}=(1+2H)\tilde{H}-H$ and {\it i) if $K<0$ and
$H>-1$, then $\tilde{\kappa}_1>0$, $\tilde{\kappa}_2>1$ and
$\tilde{H}>1$, ii) if $K>0$ and $-1<H<-1/2$, then
$0<c/2\lambda^2(1+H) <\min\{1, -H/(1+H)\}$, $\tilde{\kappa}_1<0$,
$\tilde{\kappa}_2<H/(1+H)$ and $\tilde{H}<H/(1+H)$, and iii) if
$K>0$ and $H<-1$, then $0<c/2\lambda^2(1+H) <1$,
$\tilde{\kappa}_1>(1+2H)/2(1+H)$, $\tilde{\kappa}_2>H/(1+H)$ and
$\tilde{H}>H/(1+H)$.}

The convexity of ${X_u\over|X_u|}(\log R,v)$ as a spherical curve
implies that there is a neighborhood of $O_2$ in
$\tilde{\mathcal{B}}$ which is embedded and nonnegatively curved
as a metric space. Let $\tilde{\Pi}$ be the plane parallel to
$\Pi$ and containing $\tilde{c}_1$. The curvature of $\tilde{c}_1$
is $|\tilde{\kappa}_2|/\sin\alpha$, which does not vanish. Hence
$\tilde{c}_1$ is locally convex. Using the orthogonal projection
onto $\tilde{\Pi}$, $\tilde{c}_1$ may be considered as a
$\sin\alpha$-parallel curve of $c_1$ in $\tilde{\Pi}$. Hence
$\tilde{c}_1$ is also a convex Jordan curve.

Suppose that $K<0$ and $\alpha>\pi/2$. Since $\kappa_1>0$,
$\tilde{X}_u$ is a positive multiple of $X_u$ by
\eqref{derivative}. The positivity of $\tilde{\kappa}_1$ and
$\tilde{\kappa}_2$ implies that $\tilde{\mathcal{B}}$ meets
$\tilde{\Pi}$ in constant angle $\pi-\alpha$. Suppose that $K>0$
and $\alpha <\pi/2$. If $-1<H<-1/2$, then we have $c>0$ and
$\kappa_1>-1$. Hence $\tilde{X}_u$ is a positive multiple of $X_u$
by \eqref{derivative}. The negativity of $\tilde{\kappa}_1$ and
$\tilde{\kappa}_2$ implies that $\tilde{\mathcal{B}}$ meets
$\tilde{\Pi}$ in constant angle $\alpha$. When $K>0$ and $H<-1$,
we have $c<0$ and $\kappa_1<-1$. Hence $\tilde{X}_u$ is negative
multiple of $X_u$ by \eqref{derivative}. In this case,
$\tilde{\mathcal{B}}$ lies below $\tilde{\Pi}$ and
$\tilde{\kappa}_1$ and $\tilde{\kappa}_2$ are both positive. It is
straightforward to see that $\tilde{\mathcal{B}}$ meets
$\tilde{\Pi}$ in constant angle $\alpha$.

Let $\breve{\mathcal{B}}$ be the singular surface obtained from
$\tilde{\mathcal{B}}$ by attaching the disk in $\tilde{\Pi}$
bounded by $\tilde{c}_1$ to $\tilde{\mathcal{B}}$. Since
$\tilde{\mathcal{B}}$ meets $\tilde{\Pi}$ in acute angle,
$\breve{\mathcal{B}}$ is a nonnegatively curved metric space. By
Alexandrov's generalization of Hadamard's theorem \cite{al1},
$\breve{\mathcal{B}}$ is the boundary of a convex body. Therefore
$\breve{\mathcal{B}}$ is embedded. Note again that $\tilde{H}$,
$\tilde{K}$, $\tilde{\kappa_1}$ and $\tilde{\kappa_2}$ satisfy the
statements of Lemma 2. Hence \eqref{pde} is elliptic for functions
representing $\tilde{\mathcal{B}}$ locally. We can apply
Alexandrov's moving plane argument to $\tilde{\mathcal{B}}$ using
planes perpendicular to $\tilde{\Pi}$ as in the proof of Theorem 1
to see that $\tilde{\mathcal{B}}$ is rotational. Hence
$\mathcal{B}$ is rotational and, as a result, is part of a
Delaunay surface or part of a catenoid.
\end{pf}\\

\end{document}